\documentclass[10pt]{article}
\hoffset - 20 mm
\usepackage{amsmath,amssymb,amsthm}
\usepackage{graphicx}
\usepackage{color}
\usepackage{mathtools}
\catcode`\@=11 \@addtoreset{equation}{section}

\catcode`\@=12

\newtheorem{Theorem}{Theorem}[section]
\newtheorem{Proposition}{Proposition}[section]
\newtheorem{Lemma}{Lemma}[section]
\newtheorem{Corollary}{Corollary}[section]

\newcommand{\bTheorem}[1]{
\begin{Theorem} \label{T#1} }
\newcommand{\eT}{\end{Theorem}}

\newcommand{\bProposition}[1]{
\begin{Proposition} \label{P#1}}
\newcommand{\eP}{\end{Proposition}}

\newcommand{\bLemma}[1]{
\begin{Lemma} \label{L#1} }
\newcommand{\eL}{\end{Lemma}}

\newcommand{\bCorollary}[1]{
\begin{Corollary} \label{C#1} }
\newcommand{\eC}{\end{Corollary}}

\newcommand{\bFormula}[1]{
\begin{equation} \label{#1}}
\newcommand{\eF}{\end{equation}}

\newcommand{\vr}{\varrho}

\newcommand{\vt}{\vartheta}

\newcommand{\vu}{\vc{u}}
\newcommand{\vc}[1]{{\vec #1}}
\newcommand{\Div}{{\rm div}_x}
\newcommand{\Grad}{\nabla_x}
\newcommand{\Curl}{{\rm curl}_x}

\font\F=msbm10 scaled 1000
\newcommand{\R}{\mbox{\F R}}

\newcommand{\N}{\mbox{\F N}}

\newcommand{\im}{\operatorname{Im}}
\newcommand{\cD}{{\mathcal D}}
\newcommand{\cB}{{\mathcal B}}
\newcommand{\cA}{{\mathcal A}}
\makeatletter
\DeclareRobustCommand\widecheck[1]{{\mathpalette\@widecheck{#1}}}
\def\@widecheck#1#2{%
    \setbox\z@\hbox{\m@th$#1#2$}%
    \setbox\tw@\hbox{\m@th$#1%
       \widehat{%
          \vrule\@width\z@\@height\ht\z@
          \vrule\@height\z@\@width\wd\z@}$}%
    \dp\tw@-\ht\z@
    \@tempdima\ht\z@ \advance\@tempdima2\ht\tw@ \divide\@tempdima\thr@@
    \setbox\tw@\hbox{%
       \raise\@tempdima\hbox{\scalebox{1}[-1]{\lower\@tempdima\box
\tw@}}}%
    {\ooalign{\box\tw@ \cr \box\z@}}}
\makeatother

\date{}
\makeindex
\begin{document}

\title{Global existence of a diffusion limit with damping for the compressible radiative Euler system coupled to an electromagnetic field}
\author{X. Blanc, B. Ducomet, \v S. Ne\v casov\' a}
\maketitle

 \begin{abstract}
 \noindent We study the Cauchy problem for a system of equations corresponding to a singular limit of radiative hydrodynamics, namely the 3D radiative compressible Euler system
 coupled to an electromagnetic field through the MHD approximation.
Assuming the presence of damping together with suitable smallness hypotheses for the data, we prove that this problem admits a unique global smooth solution.
\end{abstract}



\section{Introduction}

In \cite{BDN}, after the study of Buet and Despr\'es  \cite{BD} we considered a singular limit for a compressible inviscid radiative flow
where the motion of the fluid is given by the Euler system with damping for the
evolution of the density $\vr = \vr (t,x)$, the velocity
field $\vu = \vu(t,x)$, and the absolute temperature $\vt =\vt(t,x)$ as functions of the time $t$ and the Eulerian spatial
coordinate $x \in \R^3$ and 
we proved that the associated Cauchy problems admits a unique global smooth solution, provided that the data are small enough.

 In the present work we couple the previous model to an electromagnetic field through the so called magnetohydrodynamic approximation (MHD) \cite{BDu}.

Recall briefly that Maxwell's electromagnetic theory relies on Amp\`ere-Maxwell equation 
\bFormula{B1}
\partial_t \vc{D} + \vec J=\Curl \vc{H}
\eF
where $\vec D=\epsilon\vec E$ is the electric induction and  $\vec H$ is the magnetic field,
 Faraday's law 
\bFormula{b4}
\partial_t \vc{B} + \Curl \vc{E} = 0,
\eF
where  $\vc{B}=\mu \vec H$ is the magnetic induction. Here, the constant $\mu>0$ stands for the permeability of
free space. 

The two last laws are Coulomb's law
\bFormula{B2}
 \Div \vc{D} = q,
\eF
where $q$ is the electric charge density, and Gauss's law
\bFormula{B3}
 \Div \vc{B} = 0.
\eF
We first suppose that the electric current density $\vec J$ is related to the
electric field $\vc{E}$ and the macroscopic fluid velocity $\vc{u}$ via Ohm's law
\bFormula{b5}
\vc{J} = \sigma ( \vc{E} + \vc{u} \times \vc{B}),
\eF
where $\sigma$ is the electrical conductivity of the fluid.

The magnetic force acting on the fluid (Lorentz's force) $\vec f_m$ and the magnetic energy supply $E_m$ are given by
\bFormula{B4}
\vec f_m =
\vc{J} \times \vc{B},
\ \ \ E_m:=\vec J\cdot\vec E.
\eF
The MHD approximation consists in neglecting the displacement current $\partial_t \vc{D}$ (for the electric induction given by $\vec D=\epsilon\vec E$)
 in Amp\`ere-Maxwell equation and supposing that the charge $q$ is negligible, so we obtain
\bFormula{b6}
\mu \vc{J} = \Curl \vc{B},\ \mu > 0,
\eF
where, as mentioned above, the constant $\mu$ is the permeability of free space.

 Accordingly, equation (\ref{b4}) can be written \cite{C} in the form
\bFormula{b7}
\partial_t \vc{B} + \Curl( \vc{B} \times \vc{u} )
+ \Curl( \lambda \Curl \vc{B} ) = 0,
\eF
where $\lambda = (\mu \sigma)^{-1}$ is the magnetic diffusivity of the fluid.

Finally, from Faraday's law we get
\bFormula{b8}
\partial_t \left(\frac{1}{2\mu}\ |\vec B|^2\right) + \vec J\cdot\vec E = \Div \left(\frac{1}{\mu}\ \vec B\times\vec E\right).
\eF
Concerning radiation, we consider the non equilibrium diffusion regime where radiation appears through an extra equation of parabolic type
for the radiative temperature which is {\it a priori} different from the matter temperature.

More specifically the system of equations to be studied for the five unknowns $(\vr,\vu,\vt,E_r,\vec B)$ reads
\bFormula{i1bis}
\partial_t \vr + \Div (\vr \vu) = 0,
\eF
\bFormula{i2bis}
\partial_t (\vr \vu) + \Div (\vr \vu \otimes \vu) + \Grad (p+p_r)+ \frac 1 \mu \vec B\times \Curl\vec
  B +\nu\vu =0,
\eF
\bFormula{i3bis}
\partial_t \left( \vr E \right)
+ \Div \left( (\vr E+p)  \vu\right) +  \vu \cdot\Grad p_r + \frac 1 \mu \left(\vec u \times\vec B\right)\cdot \Curl\vec
  B
 = \Div \left( \kappa\Grad \vt\right)- \sigma_a\left(a\vt^4-E_r\right)+\frac \lambda \mu \left|\Curl \vec B \right|^2,
\eF
\bFormula{i4bis}
\partial_t E_r
+ \Div \left( E_r  \vu\right) + p_r\Div \vu  =  \Div \left(\frac{1}{3\sigma_s}\ \Grad E_r\right)-\sigma_a\left(E_r-a\vt^4\right),
\eF
\bFormula{i5bis}
\partial_t \vc{B} + \Curl( \vc{B} \times \vc{u} ) + \Curl\left( \lambda\ \Curl \ \vc{B}\right) = 0,\
\eF
where $\vec B$ is a divergence-free vector field, $E= \frac{1}{2} |\vu|^2 + e(\vr,\vt)$, $E_r$ is the radiative energy related
 to the temperature of radiation $T_r$ by $E_r=aT_r^4$ and $p_r$
 is the radiative pressure given by $p_r=\frac{1}{3} aT_r^4=\frac{1}{3}\ E_r$, with $a>0$.
We have also supposed for simplicity that $\mu,\sigma_a,\sigma_s,\sigma$ and $a$ are positive constants which implies in particular that 
$ \Curl\left( \frac{1}{\sigma}\ \Curl \left(\frac{1}{\mu}\ \vc{B}\right)\right)=-\frac{1}{\sigma\mu}\ \Delta\vec B$.
\vskip0.25cm
 Extending the analysis of \cite{BDN} and using stability arguments introduced by Beauchard and Zuazua in \cite{BZ}, our goal is to prove
global existence of solutions for the system  (\ref{i1bis}) - (\ref{i5bis})
 when data are sufficiently close to an equilibrium state.
\vskip0.5cm
The plan of the paper is as follows: in Section \ref{m} we state our main result (Theorem~\ref{th:existence_MHD}) then in Section \ref{MHD} we study the MHD model and prove our first result.

\section{Main result}
\label{m}

Hypotheses imposed on constitutive relations are
motivated by the general existence theory for the
Euler-Fourier system developed in \cite{Se1,Se2}. Hypotheses on transport coefficients 
are reasonable physical assumptions for the radiative part \cite{MIMI1,P}. We impose that
 pressure $p(\vr,\vt)>0$, internal energy $e(\vr,\vt)>0$ and specific entropy $s(\vr,\vt)$ are smooth functions of
 their arguments. Moreover, we impose the following monotony assumptions:
 \begin{equation}
   \label{eq:14}
\forall \vt>0, \ \forall \vr>0, \quad   \frac{\partial p}{\partial \vr}(\vr,\vt) >0, \quad \frac{\partial e}{\partial \vt}(\vr,\vt) >0.
 \end{equation}
Moreover, in our simplified setting, transport coefficients
$\kappa,\sigma_a,\sigma_s$ and the Planck's coefficient $a$ are supposed to be fixed positive numbers.
Finally the damping term with coefficient $\nu>0$ of Darcy type can be interpreted here as a diffusion of a light gas into a heavy one.
\vskip0.25cm
We are going to prove that, under the above structural assumptions on
the equation of state, system \eqref{i1bis}-\eqref{i5bis} has a global smooth solution close to any equilibrium state.

\begin{Theorem}
\label{th:existence_MHD}
  Let $\left(\overline \vr, 0, \overline\vt,\overline{E_r},\overline{\vec B}\right)$ be a constant
  state with $\overline \vr >0$, $\overline\vt >0$ and $\overline{E_r} = a{\overline \vt}^4>0$. Consider $d>7/2$.
There exists $\varepsilon>0$ such that, for any initial state
$\left(\vr_0,\vu_0,\vt_0,E_r^0,\vec B_0\right)$ satisfying
\begin{equation}
  \label{eq:donnees_petites_hors_equilibre}
  \left\|\left(\vr_0,\vu_0,\vt_0,E_r^0,\vec B_0\right) - \left(\overline \vr, 0, \overline\vt,\overline{E_r},\overline{\vec B}\right) \right\|_{H^d\left(\R^3\right)} \leq \varepsilon,
\end{equation}
there exists a unique global solution $\left(\vr,\vu,\vt,E_r,\vec B\right)$ to
\eqref{i1bis}-\eqref{i2bis}-\eqref{i3bis}-\eqref{i4bis}-\eqref{i5bis}, such that
\[
\left(\vr-\overline\vr,\vu,\vt-\overline \vt,E_r-\overline{E_r},\vec B-\overline{\vec B}\right) \in
C\left([0,+\infty);H^d\left(\R^3\right) \right)\cap C^1\left([0,+\infty);H^{d-1}\left(\R^3\right)\right).
\]
 In addition, this solution
satisfies the following energy inequality:
\begin{multline}
  \label{eq:energy_estimate_non_equilibrium}
  \left\|(\vr(t)-\overline\vr,\vu(t),\vt(t)-\overline \vt,E_r(t)-\overline{E_r},\vec B(t)-\overline{\vec B})\right\|_{H^d\left(\R^3\right)} +
  \int_0^t \left\|\Grad
      \left(\vr,\vu,\vt,E_r,\vec B\right)(s)\right\|^2_{H^{d-1}\left(\R^3\right)}ds\\
 +\int_0^t\left(
    \left\|\Grad \vt(s)\right\|^2_{H^{d}\left(\R^3\right)} + \left\|\Grad E_r(s)\right\|^2_{H^{d}\left(\R^3\right)}
 + \left\|\Grad \vec B(s)\right\|^2_{H^{d}\left(\R^3\right)}\right) ds \\
\leq C \left\|\left(\vr_0-\overline\vr,\vec u_0,\vt_0-\overline \vt,E_r^0-\overline{E_r},\vec B_0-\overline{\vec B}\right)\right\|_{H^d\left(\R^3\right)}^2,
\end{multline}
for some constant $C>0$ which does not depend on $t$.
\end{Theorem}

\section{The Euler-MHD system}
\label{MHD}

\subsection{The linearized Euler-MHD system}

Multiplying (\ref{i2bis}) by $\vec u$ and using (\ref{i1bis}) we get
\[
\partial_t \left(\frac{1}{2}\ \vr |\vu|^2\right) + \Div \left(\frac{1}{2}\ \vr |\vu|^2 \vu\right) + \Grad (p+p_r)\cdot\vu +\nu|\vu|^2 = \vec f_m\cdot\vu.
\]
Subtracting this relation from (\ref{i3bis}), using the definition $C_v=\partial_\vt e$ and the thermodynamical identity 
$\partial_\vr e=\frac{1}{\vr^2}\left(p-\vt\partial_\vt p\right)$ (Maxwell's relation),
equation (\ref{i3bis}) can be replaced by an equation for temperature
\bFormula{i3e}
\vr C_v\left(\partial_t \vt+\vu\cdot\Grad \vt\right) + \vt p_{\vt}\Div \vu-\nu\vu^2
 =  \Div \left(\kappa\Grad \vt\right)- \sigma_a\left(a\vt^4-E_r\right)+E_m-\vec f_m\cdot\vu.
\eF
Linearizing the system (\ref{i1bis})(\ref{i2bis})(\ref{i3e})(\ref{i4bis})(\ref{i5bis}) around the constant state
 $(\overline{\vr},0,\overline{\vt},\overline{E}_r,{\overline{\vec B}})$
with the compatibility condition $\overline{E}_r=a\overline{\vt}^4$
 and putting $\vr=r+\overline{\vr}$, $\vt=T+\overline{\vt}$, $E_r=e_r+\overline{E}_r$ and $\vec B=\vec b+\overline{\vec B}$ we get
\bFormula{i1bislin}
\partial_t r +\overline{\vr}\ \Div \vu = 0,
\eF
\bFormula{i2bislin}
\partial_t  \vu
+ \frac{\overline{p}_{\vr}}{\overline{\vr}} \Grad r
+\frac{\overline{p}_{\vt}}{\overline{\vr}}\Grad T
+\frac{1}{3\overline{\vr}}\Grad e_r +\ \frac{1}{\mu}\overline{\vec B}\times \left(\Curl \vec b\right)+\nu\vu = 0,
\eF
\bFormula{i3bislin}
\partial_t T
+ \frac{\overline{\vt}\overline{p}_{\vt}}{\overline{\vr}\overline{C}_v}\ \Div \vu
=  \Div \left(
\frac{\kappa}{\overline{\vr}\overline{C}_v} \Grad T\right)
-\frac{\sigma_a}{\overline{\vr}\overline{C}_v}\left(4a\overline{\vt}^3T-e_r\right),
\eF
\bFormula{i4bislin}
\partial_t e_r
+ \frac{4}{3} \overline{E}_r \Div \vu  =  \Div \left(\frac{1}{3\sigma_s}\ \Grad e_r\right)-\sigma_a\left(e_r-4a\overline{\vt}^3T\right),
\eF
\bFormula{I4bislin}
\partial_t \vec b +\overline{\vec B}\Div \vu-(\overline{\vec B}\cdot\Grad)\vu=  \lambda\ \Delta \vec b.
\eF
Using the vector notation $U:=
\left(
 \begin{array}{c}
{\displaystyle r}\\\
{\displaystyle  u_1}\\\
{\displaystyle  u_2}\\\
{\displaystyle  u_3}\\\
{\displaystyle  T}\\\
{\displaystyle  e_r}\\\
{\displaystyle  b_1}\\\
{\displaystyle  b_2}\\\
{\displaystyle  b_3}\\\
\end{array}
\right)$, the linearized system (\ref{i1bislin}) - (\ref{I4bislin}) rewrites
\begin{equation}
{\displaystyle \partial_t U+\sum_{j=1}^3 \cA_j\partial_j U =\cD\Delta U-\cB U,}
\label{hrad}
\end{equation}
with
\[
\cA_1=
\left(
\begin{array}{ccccccccc}
0  &  \overline{\vr}  &  0  &  0 & 0 & 0 & 0 & 0 & 0\\\
\alpha' & 0 & 0 & 0 & \beta' & \beta'' & 0 & \overline{B}_2/\mu & \overline{B}_3/\mu\\\
0 & 0 & 0 & 0 & 0 & 0 & 0 & -\overline{B}_1/\mu & 0\\\
0 & 0 & 0 & 0 & 0 & 0 & 0 & 0 & -\overline{B}_1/\mu\\\
0 & \gamma' & 0 & 0 & 0 & 0 & 0 & 0 & 0\\\
0 & \gamma'' & 0 & 0 & 0 & 0 & 0 & 0 & 0\\\
0 & 0 & 0 & 0 & 0 & 0 & 0 & 0 & 0\\\
0 & \overline{B}_2 & -\overline{B}_1 & 0 & 0 & 0 & 0 & 0 & 0\\\
0 & \overline{B}_3 & 0 & -\overline{B}_1 & 0 & 0 & 0 & 0 & 0
\end{array}
\right),
\]
\medskip
\[
\cA_2=
\left(
\begin{array}{ccccccccc}
0  & 0 &  \overline{\vr}  &  0  &  0 & 0 & 0 & 0 & 0 \\\
0 & 0 & 0 & 0 & 0 & 0 & -\overline{B}_2/\mu & 0 & 0\\\
\alpha' & 0 & 0 & 0 & \beta' & \beta'' & \overline{B}_1/\mu & 0 & \overline{B}_3/\mu\\\
0 & 0 & 0 & 0 & 0 & 0 & 0 & 0 & -\overline{B}_2/\mu\\\
0 & 0 & \gamma' & 0 & 0 & 0 & 0 & 0 & 0 \\\
0 & 0 & \gamma'' & 0 & 0 & 0 & 0 & 0 & 0\\\
0 & -\overline{B}_2 & \overline{B}_1 & 0 & 0 & 0 & 0 & 0 & 0\\\
0 & 0 & 0 & 0 & 0 & 0 & 0 & 0 & 0\\\
0 & 0 & \overline{B}_3 & -\overline{B}_2 & 0 & 0 & 0 & 0 & 0
\end{array}
\right)
\ \]
\medskip
\[
\cA_3=
\left(
\begin{array}{ccccccccc}
0 & 0 & 0 & \overline{\vr}  &  0  &  0 & 0 & 0 & 0\\\
0 & 0 & 0 & 0 & 0 & 0 & -\overline{B}_3/\mu & 0 & 0\\\
0 & 0 & 0 & 0 & 0 & 0 & 0 & -\overline{B}_3/\mu & 0\\\
\alpha' & 0 & 0 & 0 & \beta' & \beta'' & \overline{B}_1/\mu & \overline{B}_2/\mu & 0\\\
0 & 0 & 0 & \gamma' & 0 & 0 & 0 & 0 & 0\\\
0 & 0 & 0 & \gamma'' & 0 & 0 & 0 & 0 & 0\\\
0 & -\overline{B}_3 & 0 & \overline{B}_1 & 0 & 0 & 0 & 0 & 0\\\
0 & 0 & -\overline{B}_3 & \overline{B}_2 & 0 & 0 & 0 & 0 & 0\\\
0 & 0 & 0 & 0 & 0 & 0 & 0 & 0 & 0
\end{array}
\right),
\]
and
\[
\cD=
\left(
\begin{array}{ccccccccc}
0 & 0 & 0 & 0 & 0 & 0 & 0 & 0 & 0\\\
0 & 0 & 0 & 0 & 0 & 0 & 0 & 0 & 0\\\
0 & 0 & 0 & 0 & 0 & 0 & 0 & 0 & 0\\\
0 & 0 & 0 & 0 & 0 & 0 & 0 & 0 & 0\\\
0 & 0 & 0 & 0 & \delta' & 0 & 0 & 0 & 0\\\
0 & 0 & 0 & 0 & 0 & \delta'' & 0 & 0 & 0\\\
0 & 0 & 0 & 0 & 0 & 0 & \lambda & 0 & 0\\\
0 & 0 & 0 & 0 & 0 & 0 & 0 & \lambda & 0\\\
0 & 0 & 0 & 0 & 0 & 0 & 0 & 0 & \lambda
\end{array}
\right)
\ \
\cB:=
\left(
\begin{array}{ccccccccc}
0 & 0 & 0 & 0 & 0 & 0 & 0 & 0 & 0\\\
0 & \nu  & 0 & 0 & 0 & 0 & 0 & 0 & 0\\\
0 & 0 & \nu  & 0 & 0 & 0 & 0 & 0 & 0\\\
0 & 0 & 0 & \nu  & 0 & 0 & 0 & 0 & 0\\\
0 & 0 & 0 & 0 & \zeta & -\eta & 0 & 0 & 0\\\
0 & 0 & 0 & 0 & -\pi & \sigma_a & 0 & 0 & 0\\\
0 & 0 & 0 & 0 & 0 & 0 & 0 & 0 & 0\\\
0 & 0 & 0 & 0 & 0 & 0 & 0 & 0 & 0\\\
0 & 0 & 0 & 0 & 0 & 0 & 0 & 0 & 0
\end{array}
\right),
\]
where
\[
\alpha'= \frac{\overline{p}_{\vr}}{\overline{\vr}},\quad
\beta'=\frac{\overline{p}_{\vt}}{\overline{\vr}},\quad
\beta''=\frac{1}{3\overline{\vr}},\quad
\gamma'=\frac{\overline{p}_{\vt}}{\overline{C}_v} ,\quad
\delta'= \frac{\kappa}{\overline{\vr}\overline{C}_v},\]
\[
\gamma'' = \frac43 \overline E_r,\quad
\delta''=\frac{1}{3\sigma_s},\quad
\zeta= \frac{4a\sigma_a\overline{\vt}^3}{\overline{\vr}\overline{C}_v},\quad
\eta=\frac{\sigma_a}{\overline{\vr}\overline{C}_v},\quad
\pi=4a\sigma_a\overline{\vt}^3.
\]
In order to apply the Kreiss theory we have to put the system (\ref{hrad}) in a symmetric form \cite{BGS}.
For that purpose it is sufficient to consider a diagonal symmetrizer
\begin{equation}
\widetilde \cA_0=
\left(
\begin{array}{ccccccccc}
\mu\frac{\alpha'}{\overline\vr} & 0 & 0 & 0 & 0 & 0 & 0 & 0 & 0\\\
0 & \mu & 0 & 0 & 0 & 0 & 0 & 0 & 0\\\
0 & 0 & \mu & 0 & 0 & 0 & 0 & 0 & 0\\\
0 & 0 & 0 & \mu & 0 & 0 & 0 & 0 & 0\\\
0 & 0 & 0 & 0 & \mu\frac{\beta'}{\gamma'} & 0 & 0 & 0 & 0\\\
0 & 0 & 0 & 0 & 0 & \mu\frac{\beta''}{\gamma''} & 0 & 0 & 0\\\
0 & 0 & 0 & 0 & 0 & 0 & 1 & 0 & 0\\\
0 & 0 & 0 & 0 & 0 & 0 & 0 & 1 & 0\\\
0 & 0 & 0 & 0 & 0 & 0 & 0 & 0 & 1
\end{array}
\right).
\label{a0}
\end{equation}

Multiplying the first equation (\ref{hrad}) by $\widetilde \cA_0$ on the left, we get
\begin{equation}
{\displaystyle \widetilde \cA_0\partial_t U+\sum_{j=1}^3 \widetilde \cA_j\partial_j U =\widetilde \cD\Delta
  U-\widetilde \cB U,}
\label{hsym}
\end{equation}
where the matrices $\widetilde \cA_j=\widetilde \cA_0 \cA_j$ are symmetric,
for all  $j=1,2,3$. More specifically,
\[
\widetilde \cA_1=
\left(
\begin{array}{ccccccccc}
0  &  \mu\alpha'  &  0  &  0 & 0 & 0 & 0 & 0 & 0\\\
\mu\alpha' & 0 & 0 & 0 & \mu\beta' & \mu\beta'' & 0 & \overline{B}_2 & \overline{B}_3\\\
0 & 0 & 0 & 0 & 0 & 0 & 0 & -\overline{B}_1 & 0\\\
0 & 0 & 0 & 0 & 0 & 0 & 0 & 0 & -\overline{B}_1\\\
0 & \mu\beta' & 0 & 0 & 0 & 0 & 0 & 0 & 0\\\
0 & \mu\beta'' & 0 & 0 & 0 & 0 & 0 & 0 & 0\\\
0 & 0 & 0 & 0 & 0 & 0 & 0 & 0 & 0\\\
0 & \overline{B}_2 & -\overline{B}_1 & 0 & 0 & 0 & 0 & 0 & 0\\\
0 & \overline{B}_3 & 0 & -\overline{B}_1 & 0 & 0 & 0 & 0 & 0
\end{array}
\right),
\]
\smallskip
\[
\widetilde \cA_2=
\left(
\begin{array}{ccccccccc}
0  & 0 &  \mu\alpha'  &  0  &  0 & 0 & 0 & 0 & 0 \\\
0 & 0 & 0 & 0 & 0 & 0 & -\overline{B}_2 & 0 & 0\\\
\mu\alpha' & 0 & 0 & 0 & \mu\beta' & \mu\beta'' & \overline{B}_1 & 0 & \overline{B}_3\\\
0 & 0 & 0 & 0 & 0 & 0 & 0 & 0 & -\overline{B}_2\\\
0 & 0 & \mu\beta' & 0 & 0 & 0 & 0 & 0 & 0 \\\
0 & 0 & \mu\beta'' & 0 & 0 & 0 & 0 & 0 & 0\\\
0 & -\overline{B}_2 & \overline{B}_1 & 0 & 0 & 0 & 0 & 0 & 0\\\
0 & 0 & 0 & 0 & 0 & 0 & 0 & 0 & 0\\\
0 & 0 & \overline{B}_3 & -\overline{B}_2 & 0 & 0 & 0 & 0 & 0
\end{array}
\right)
,\]
\smallskip
\[
\widetilde \cA_3=
\left(
\begin{array}{ccccccccc}
0 & 0 & 0 & \mu\alpha'  &  0  &  0 & 0 & 0 & 0\\\
0 & 0 & 0 & 0 & 0 & 0 & -\overline{B}_3 & 0 & 0\\\
0 & 0 & 0 & 0 & 0 & 0 & 0 & -\overline{B}_3 & 0\\\
\mu\alpha' & 0 & 0 & 0 & \mu\beta' & \mu\beta'' & \overline{B}_1 & \overline{B}_2 & 0\\\
0 & 0 & 0 & \mu\beta' & 0 & 0 & 0 & 0 & 0\\\
0 & 0 & 0 & \mu\beta'' & 0 & 0 & 0 & 0 & 0\\\
0 & -\overline{B}_3 & 0 & \overline{B}_1 & 0 & 0 & 0 & 0 & 0\\\
0 & 0 & -\overline{B}_3 & \overline{B}_2 & 0 & 0 & 0 & 0 & 0\\\
0 & 0 & 0 & 0 & 0 & 0 & 0 & 0 & 0
\end{array}
\right).
\]
The hyperbolic part of system (\ref{hsym}) is now symmetric while its symmetric dissipative part is given by
 \begin{equation}
   \label{eq:def_BD}
\widetilde \cD=
\left(
\begin{array}{ccccccccc}
0 & 0 & 0 & 0 & 0 & 0 & 0 & 0 & 0\\\
0 & 0 & 0 & 0 & 0 & 0 & 0 & 0 & 0\\\
0 & 0 & 0 & 0 & 0 & 0 & 0 & 0 & 0\\\
0 & 0 & 0 & 0 & 0 & 0 & 0 & 0 & 0\\\
0 & 0 & 0 & 0 & \frac{\mu\beta'\delta'}{\gamma'} & 0 & 0 & 0 & 0\\\
0 & 0 & 0 & 0 & 0 & \frac{\mu\beta''\delta''}{\gamma''} & 0 & 0 & 0\\\
0 & 0 & 0 & 0 & 0 & 0 & \lambda & 0 & 0\\\
0 & 0 & 0 & 0 & 0 & 0 & 0 & \lambda & 0\\\
0 & 0 & 0 & 0 & 0 & 0 & 0 & 0 & \lambda
\end{array}
\right)
\ \widetilde \cB =
\left(
\begin{array}{ccccccccc}
0 & 0 & 0 & 0 & 0 & 0 & 0 & 0 & 0\\\
0 &\mu \nu  & 0 & 0 & 0 & 0 & 0 & 0 & 0\\\
0 & 0 & \mu\nu  & 0 & 0 & 0 & 0 & 0 & 0\\\
0 & 0 & 0 & \mu\nu  & 0 & 0 & 0 & 0 & 0\\\
0 & 0 & 0 & 0 & \frac{\mu\beta'\zeta}{\gamma'} & -\frac{\mu\beta'\eta}{\gamma'} & 0 & 0 & 0\\\
0 & 0 & 0 & 0 & -\frac{\mu\beta''\pi}{\gamma''} & \frac{\mu\sigma_a\beta''}{\gamma''} & 0 & 0 & 0\\\
0 & 0 & 0 & 0 & 0 & 0 & 0 & 0 & 0\\\
0 & 0 & 0 & 0 & 0 & 0 & 0 & 0 & 0\\\
0 & 0 & 0 & 0 & 0 & 0 & 0 & 0 & 0
\end{array}
\right),
 \end{equation}
where one checks the positiveness condition of $\widetilde \cB $
\[
 ^t X\widetilde \cB X\geq 0,\ \ \mbox{for any vector}\  X\in\R^9.
\]
Applying the Fourier transform in $x$ to (\ref{hsym}) we get
\begin{equation}
{\displaystyle \widetilde \cA_0\partial_t \widehat U+i\sum_{j=1}^3 \xi_j\widetilde \cA_j\widehat U = -|\xi|^2\widetilde \cD\widehat U -\widetilde \cB \widehat U,}
\label{hsym1}
\end{equation}
or
\begin{equation}
{\displaystyle \widetilde \cA_0\partial_t \widehat U=E(\xi)\widehat U,}
\label{hsym2}
\end{equation}
with
\[
E(\xi)=-\cB (\xi)-i\cA(\xi),
\]
where 
\[
\cA(\xi)=\sum_{j=1}^3 \xi_j\widetilde \cA_j=
\]
\begin{equation}\label{eq:3}
\left(
\begin{array}{ccccccccc}
0 & \mu\alpha'\xi_1 & \mu\alpha'\xi_2 & \mu\alpha'\xi_3 & 0 & 0 & 0 & 0 & 0\\\
\mu\alpha'\xi_1 & 0 & 0 & 0 & \mu\beta'\xi_1 & \mu\beta''\xi_1 & -\overline{B}_2\xi_2-\overline{B}_3\xi_3 & \overline{B}_2\xi_1 & \overline{B}_3\xi_1\\\
\mu\alpha'\xi_2 & 0 & 0 & 0 & \mu\beta'\xi_2 & \mu\beta''\xi_2 & \overline{B}_1\xi_2 & -\overline{B}_1\xi_1 -\overline{B}_3\xi_3 & \overline{B}_3\xi_2\\\
\mu\alpha'\xi_3 & 0 & 0 & 0 & \mu\beta'\xi_3 & \mu\beta''\xi_3 & \overline{B}_1\xi_3 & \overline{B}_2\xi_3 &  -\overline{B}_1\xi_1- \overline{B}_2\xi_2\\\
0 & \mu\beta' \xi_1 & \mu\beta'\xi_2 & \mu\beta' \xi_3 & 0 & 0 & 0 & 0 & 0\\\
0 & \mu\beta'' \xi_1 & \mu\beta''\xi_2 & \mu\beta'' \xi_3 & 0 & 0 & 0 & 0 & 0\\\
0 & -\overline{B}_2\xi_2-\overline{B}_3\xi_3 & \overline{B}_1\xi_2 & \overline{B}_1\xi_3 & 0 & 0 & 0 & 0 & 0\\\
0 & \overline{B}_2\xi_1 & -\overline{B}_1\xi_1- \overline{B}_3\xi_3 & \overline{B}_2\xi_3 & 0 & 0 & 0 & 0 & 0\\\
0 & \overline{B}_3\xi_1 & \overline{B}_3\xi_2 & -\overline{B}_1\xi_1- \overline{B}_2\xi_2 & 0 & 0 & 0 & 0 & 0
\end{array}
\right)
\end{equation}

and
\begin{equation}\label{eq:3bis}
\cB (\xi):=\widetilde \cB +|\xi|^2\widetilde \cD =
\left(
\begin{array}{ccccccccc}
0 & 0 & 0 & 0 & 0 & 0 & 0 & 0 & 0\\\
0 & \nu  & 0 & 0 & 0 & 0 & 0 & 0 & 0\\\
0 & 0 &  \nu & 0 & 0 & 0 & 0 & 0 & 0\\\
0 & 0 & 0 &  \nu & 0 & 0 & 0 & 0 & 0\\\
0 & 0 & 0 & 0 & \frac{\mu\beta'\zeta}{\gamma'}+\frac{\mu\beta'\delta'}{\gamma'}|\xi|^2  & -\frac{\mu\beta'\eta}{\gamma'} & 0 & 0 & 0\\\
 0 & 0 & 0 & 0 & -\frac{\mu\beta''\pi}{\gamma''} & \frac{\mu\beta''\sigma_a}{\gamma''}+\frac{\mu\beta''\delta''}{\gamma''}|\xi|^2 & 0 & 0 & 0\\\
0 & 0 & 0 & 0 & 0 & 0 & \lambda|\xi|^2 & 0 & 0\\\
0 & 0 & 0 & 0 & 0 & 0 & 0 & \lambda|\xi|^2 & 0\\\
0 & 0 & 0 & 0 & 0 & 0 & 0 & 0 & \lambda|\xi|^2
\end{array}
\right).
\end{equation}
Solving this equation with initial condition $\widehat U_0(\xi)$ we get
\begin{equation}
{\displaystyle  \widehat U(t,\xi)=\exp\left[t\widetilde \cA_0^{-1}E(\xi)\right]\widehat U_0(\xi).}
\label{sol}
\end{equation}
In the strictly hyperbolic case $\widetilde \cD=0$, under the Kalman rank
condition \cite{KLA} for the pair $(\cA(\xi),\cB )$, it can be proved \cite{BZ}
that
\[
\exists C>0,\ \ \lambda(\xi)>0\ \ :\ \ \exp\left[t\widetilde \cA_0^{-1}E(\xi)\right]\leq C e^{-\lambda(\xi)t}.
\]
Observing the partially parabolic character of the system, one can expect a similar result when $\widetilde \cD\neq 0$ with a parabolic smoothing effect at
low frequencies and an extra damping in the high frequency regime.

Taking benefit of the damping, we can use the Shizuta-Kawashima condition (SK) \cite{SK} which applies to the previous system.
Following the arguments of Beauchard and Zuazua \cite{BZ}, we have
\newtheorem{l1}{Lemma}
\begin{l1} For any $\xi \in S^2,$ a necessary and sufficient condition for the matrices $\cB (\xi)$ and $\cA(\xi)$ defined by (\ref{eq:3}) and (\ref{eq:3bis}) satisfy the Shizuta-Kawashima condition (SK):
  \begin{equation}
    \label{eq:kawashima_shizuta}
\left\{\mbox{eigenvectors of}\ \left(\widetilde \cA_0\right)^{-1}\cA(\xi) \right\}\cap \ker\ \cB (\xi)=\{0\},
  \end{equation}
is that $\nu>0$.
\label{l1}
\end{l1}

 {\bf Proof:} 
\begin{enumerate}
\item  If $\nu\neq 0$. One checks that $\ker\ \cB (\xi)$ is the
1-dimensional subspace spanned by the vector $(1,0,0,0,0,0,0,0,0)$. Therefore,
if $X\in \ker\ \cB (\xi)\setminus\{0\}$ is an eigenvector of $\left(\widetilde
  \cA_0\right)^{-1}\cA(\xi)$, we have $X = (x_1,0,0,0,0,0,0,0,0)$, $x_1\neq 0$, and
\[
\cA(\xi) X= \lambda \widetilde \cA_0 X,
\]
for some $\lambda\in\R$. According to the values of $\widetilde \cA_0$ and $\cA(\xi)$,
this implies that $\lambda = 0$, $\xi_1 = \xi_2 = \xi_3 = 0$, which is
in contradiction with the hypotesis $\xi\in S^2$.

\item  If $\nu= 0$. One checks that $\ker\ \cB (\xi)$ is the 4-dimensional subspace spanned by the vectors $(x_1,x_2,x_3,x_4,0,0,0,0,0)$.

 Let us denote by $(\lambda,X)$ an eigenpair of $\cA(\xi)$, with non zero eigenvector $X\in  \ker \cB (\xi)$.
$X$ satisfies the system
\[
\mu\alpha'\xi_1 x_2+ \mu\alpha'\xi_2x_3 + \mu\alpha'\xi_3x_4=\lambda x_1,
\]
\[
\mu\alpha'\xi_1 x_1=\lambda x_2,
\]
\[
\mu\alpha'\xi_2 x_1=\lambda x_3,
\]
\[
\mu\alpha'\xi_3 x_1=\lambda x_4,
\]
\[
\mu\beta' \xi_1 x_2+ \mu\beta'\xi_2 x_3+ \mu\beta' \xi_3 x_4=0,
\]
\[
\mu\beta'' \xi_1 x_2+ \mu\beta''\xi_2 x_3+ \mu\beta'' \xi_3 x_4=0,
\]
\[
 -(\overline{B}_2\xi_2+\overline{B}_3\xi_3)x_2+  \overline{B}_1\xi_2 x_3+ \overline{B}_1\xi_3 x_4=0,
\]
\[
\overline{B}_2\xi_1 x_2 -(\overline{B}_1\xi_1+ \overline{B}_3\xi_3)x_3+ \overline{B}_2\xi_3 x_4=0,
\]
\[
 \overline{B}_3\xi_1 x_2+ \overline{B}_3\xi_2 x_3 -(\overline{B}_1\xi_1+ \overline{B}_2\xi_2) x_4=0.
\]
Denoting $\vec B=(\overline{B}_1,\overline{B}_2,\overline{B}_3)$, $\vec \xi=(\xi_1,\xi_2,\xi_3)$, and $\vec x=(x_2,x_3,x_4)$, the system rewrites
\[
\mu\alpha'\ \vec x\cdot \vec\xi=\lambda x_1,
\]
\[
\mu\alpha' x_1\vec \xi=\lambda \vec x,
\]
\[
\mu\beta'\ \vec x\cdot \vec\xi=0,
\]
\[
\mu\beta''\ \vec x\cdot \vec\xi=0,
\]
\[ 
-(\vec B\cdot \vec\xi)\vec x+(\vec \xi\cdot\vec x)\vec B=0.
\]
In particular, this implies $\vec x \cdot\vec \xi = 0$, which in turn implies $\lambda x_1 = 0$. As a consequence,
we have
\[
\lambda x_1=0,
\]
\[
\vec x\cdot \vec\xi=0,
\]
\[
\mu\alpha' x_1\vec \xi=\lambda \vec x,
\]
\[ 
(\vec \xi\cdot\vec B)\vec x=0.
\]
Choosing $\lambda=0$, we see that any $\vec x\in \vec \xi^{\perp}$ for $\xi\in \vec B^{\perp}$ gives a nontrivial eigenpair $(\lambda,X)$ with $\lambda=0$ and 
$X=(0,\vec x,0,0,0,0,0)$. Then the SK condition is not satisfied.

\end{enumerate}

As in the equilibrium case \eqref{eq:kawashima_shizuta} is equivalent to the
existence of a compensating matrix:

\begin{Proposition}
  \label{pr:compensating_matrix}
For any $\xi \in S^2,$ the matrices $\widetilde \cA_0$, $\cB (\xi)$ and $\cA(\xi)$ being defined by (\ref{a0}), (\ref{eq:3}) and (\ref{eq:3bis}),
there exists a matrix-valued function
\begin{displaymath}
  \begin{array}{rrcl}
    K : & S^2 & \longrightarrow & \R^{6\times 6} \\
    & \omega & \longmapsto & K(\omega)
  \end{array}
\end{displaymath}
such that
\begin{enumerate}
\item $\omega\mapsto K(\omega)$ is a $C^\infty$ function, and satisfies
  $K(-\omega) = - K(\omega)$ for any $\omega \in S^2$.
\item $K(\omega) \tilde \cA_0$ is a skew-symmetric matrix for any
  $\omega\in S^2$.
\item Denoting by $[A] = \frac12 \left(A+ A^T\right)$ the symmetric part
  of $A$, the matrix $[K(\omega)\cA(\omega)] + \cB (\omega)$ is symmetric
  positive definite for any $\omega\in S^2$.
\end{enumerate}
\end{Proposition}

\subsection{Entropy properties}

Adding equations \eqref{b8}, (\ref{i3bis}) and (\ref{i4bis}) we get
\begin{multline}
  \label{iener}
\partial_t \left( \frac{1}{2}\vr \left|\vu\right|^2+\vr e +E_r+\frac{1}{2\mu}\ |\vec B|^2\right)
+ \Div \left( (\vr E+E_r)  \vu + (p+ p_r)\vu+\frac{1}{\mu}\ \vec E\times\vec B\right)\\
 = \Div \left( \kappa\Grad \vt\right)+\Div \left(\frac{1}{3\sigma_s}\ \Grad E_r\right).  
\end{multline}
Introducing the entropy $s$ of the fluid by the Gibbs law $\vt ds=de+pd\left(\frac{1}{\vr}\right)$ and denoting by $S_r:=\frac{4}{3}aT_r^3$ the radiative entropy,
 equation (\ref{i4bis}) rewrites
\[
\partial_t S_r
+ \Div \left( S_r  \vu\right)  = \frac{1}{T_r} \Div \left(\frac{1}{3\sigma_s}\ \Grad E_r\right)-\sigma_a\frac{E_r-a\vt^4}{T_r},
\]
or
\bFormula{iER}
\partial_t S_r
+ \Div \left( S_r  \vu\right)  =  \Div \left(\frac{1}{3\sigma_sT_r}\
  \Grad E_r\right)+\frac{4a}{3\sigma_s}T_r|\Grad T_r|^2-\sigma_a\frac{E_r-a\vt^4}{T_r}.
\eF
Replacing equation (\ref{i3bis}) by the internal energy equation
\begin{equation}
\partial_t (\vr e) + \Div (\vr e \vu) + p\Div \vu - \nu \left|\vu\right|^2
= \Div \left( \kappa\Grad \vt\right)- \sigma_a\left(a\vt^4-E_r\right)+\frac{1}{\sigma\mu^2}\ |\Curl\vec B|^2,
\label{intern}
\end{equation}
and dividing it by $\vt$, we may write an entropy equation for matter
\begin{equation}
\partial_t (\vr s) + \Div (\vr s \vu) -\frac\nu\vt\left|\vu\right|^2 = \Div \left( \frac{\kappa\Grad\vt}{\vt} \right)
 +\frac{\kappa|\Grad\vt|^2}{\vt^2} - \sigma_a\frac{a\vt^4-E_r}{\vt}+\frac{1}{\sigma\mu^2\vt}\ |\Curl\vec B|^2.
\label{entro}
\end{equation}
So adding (\ref{entro}) and (\ref{iER}) we obtain
\[
\partial_t \left( \vr  s +S_r\right)
+ \Div \left( (\vr  s  +S_r)\vu\right)
-\Div \left(   \frac{\kappa\Grad\vt}{\vt} +\frac{1}{3\sigma_sT_r}\ \Grad E_r   \right)
\]
\bFormula{ient}
=\frac{\kappa|\Grad\vt|^2}{\vt^2}+\frac{4a}{3\sigma_s}T_r|\Grad E_r|^2+
\frac{a\sigma_a}{\vt
  T_r}\left(\vt-T_r\right)^2\left(\vt+T_r\right)\left(\vt^2+T_r^2\right)+\frac{1}{\sigma\mu^2\vt}\ |\Curl\vec B|^2
+ \frac\nu\vt \left|\vu\right|^2.
\eF
Introducing the Helmholtz functions $H_{\overline{\vt}}(\vr,\vt):=\vr\left(e-\overline{\vt}s\right)$ and $H_{r,\overline{\vt}}(T_r):=E_r-\overline{\vt}S_r$, we check that
the quantities $H_{\overline{\vt}}(\vr,\vt)-(\vr-\overline{\vr})\partial_{\vr}H_{\overline{\vt}}(\overline{\vr},\overline{\vt})-H_{\overline{\vt}}(\overline{\vr},\overline{\vt})$ and
$H_{r,\overline{\vt}}(T_r)-H_{r,\overline{\vt}}(\overline{T}_r)$ are non-negative and strictly coercive functions reaching zero minima at the equilibrium state $(\overline{\vr},\overline{\vt},\overline{E}_r)$.

\newtheorem{l2}[l1]{Lemma}
\begin{l2}
Let $\overline{\vr}$ and $\overline{\vt} = \overline{T}_r$ be given positive constants.
Let ${\mathcal O}_1$ and ${\mathcal O}_2$ be the sets defined by

\bFormula{es1}
{\mathcal O}_1:=\left\{
(\vr,\vt)\in\R^2\ :\ \frac{\overline{\vr}}{2}<\vr<2\overline{\vr},\ \frac{\overline{\vt}}{2}<\vt<2\overline{\vt},
\right\}.
\eF
\bFormula{es2}
{\mathcal O}_2:=\left\{
 T_r\in\R\ :\ \frac{\overline{T}_r}{2}<T_r<2\overline{T}_r,
\right\}.
\eF
There exist positive constants $C_{1,2}(\overline{\vr},\overline{\vt})$ and $C_{3,4}(\overline{T}_r)$ such that
\begin{enumerate}
\item\label{item:1}
\[
C_1\left(|\vr-\overline{\vr}|^2+|\vt-\overline{\vt}|^2\right)
\leq
H_{\overline{\vt}}(\vr,\vt)
-(\vr-\overline{\vr})\partial_{\vr}H_{\overline{\vt}}(\overline{\vr},\overline{\vt})
-H_{\overline{\vt}}(\overline{\vr},\overline{\vt})
\]
\bFormula{ess1}
\leq
C_2\left(|\vr-\overline{\vr}|^2+|\vt-\overline{\vt}|^2\right),
\eF
for all $(\vr,\vt)\in{\mathcal O}_1$,
\item\label{item:2}
\bFormula{ess2}
C_3 |T_r-\overline{T}_r|^2
\leq
 H_{r,\overline{\vt}}(T_r)-H_{r,\overline{\vt}}(\overline{T}_r)
\leq
C_4|T_r-\overline{T}_r|^2,
\eF
for all $T_r\in{\mathcal O}_2$.
\end{enumerate}
\label{l2}
\end{l2}

\noindent {\bf Proof:}
\begin{enumerate}
\item Point~\ref{item:1} is proved in \cite{FEINOV} and we only sketch the proof for convenience.

We have the decomposition
$$\vr\mapsto
H_{\overline{\vt}}(\vr,\vt)
-(\vr-\overline{\vr})\partial_{\vr}H_{\overline{\vt}}(\overline{\vr},\overline{\vt})
-H_{\overline{\vt}}(\overline{\vr},\overline{\vt})
={\mathcal F}(\vr)+{\mathcal G}(\vr),$$
where
${\mathcal F}(\vr)=H_{\overline{\vt}}(\vr,\overline{\vt})
-(\vr-\overline{\vr})\partial_{\vr}H_{\overline{\vt}}(\overline{\vr},\overline{\vt})
-H_{\overline{\vt}}(\overline{\vr},\overline{\vt})$ and
 $ {\mathcal G}(\vr)=H_{\overline{\vt}}(\vr,\vt)-H_{\overline{\vt}}(\vr,\overline{\vt})$. Using Gibbs law $\vt ds =
 de + pd\left(\frac 1 \vr\right)$, one easily proves that $\partial^2_\vr H_{\overline \vt}(\vr,\overline \vt)=
 \frac 1 {\vr \overline\vt}\partial_\vr p(\vr,\overline\vt),$ which is positive according to \eqref{eq:14}. Hence,
 ${\mathcal F}$ is strictly convex and reaches a zero minimum at $\overline{\vr}$. Turning to $\mathcal{G}$, we
 have, still using Gibbs law, $\partial_\vt H_{\overline \vt}(\vr,\overline\vt) =
 \vr\frac{\vt-\overline\vt}\vt \partial_\vt e(\vr,\overline\vt)$.
 Thus, using \eqref{eq:14} again, we infer that $ {\mathcal G}$ is strictly
decreasing for $\vt<\overline{\vt}$ and strictly
increasing for $\vt>\overline{\vt}$.
Computing the derivatives of $H_{\overline{\vt}}$ leads directly to the estimate (\ref{ess1}).

\item Point~\ref{item:2} follows from the properties of the function
$x\mapsto H_{r,\overline{\vt}}(x)-H_{r,\overline{\vt}}\left(\overline T_r\right)=ax^3(x-\frac{4}{3}\overline{\vt})+\frac{a}{3}\overline{\vt}^4$\hfill$\Box$
\end{enumerate}
From this simple result, we can obtain an identity leading to energy estimates.
In fact, multiplying (\ref{ient}) by $\overline\vt$, subtracting the result to (\ref{iener}) and using the conservation of mass, we get
\[
\partial_t \left( \frac{1}{2}\vr \left|\vu\right|^2
+H_{\overline{\vt}}(\vr,\vt)-(\vr-\overline{\vr})\partial_{\vr}H_{\overline{\vt}}(\overline{\vr},\overline{\vt})-H_{\overline{\vt}}(\overline{\vr},\overline{\vt})
+H_{r,\overline{\vt}}(T_r) - H_{r,\overline\vt}\left(\overline T_r\right) +\frac 1 {2\mu}|\vc{B}|^2
 \right)
\]
\[
+ \Div \left( \left(\vr (E - \overline e) +E_r\right)  \vu + (p+ p_r)\vu - \overline{\vt} (\vr (s-\overline s)+S_r)  \vu  +\frac 1 \mu \vc{E}\times \vc{B}\right)
\]
\[
 = \Div \left( \kappa\Grad \vt+\frac{1}{3\sigma_s}\ \Grad E_r\right)
-\overline{\vt}\Div \left(   \frac{\kappa\Grad\vt}{\vt} +\frac{1}{3\sigma_sT_r}\ \Grad E_r   \right)
\]
\bFormula{ienerbis}
-\overline{\vt}\frac{\kappa|\Grad\vt|^2}{\vt^2}
-\overline{\vt}\frac{4a}{3\sigma_s}T_r|\Grad E_r|^2
- \overline{\vt}\frac{a\sigma_a}{\vt
  T_r}\left(\vt-T_r\right)^2\left(\vt+T_r\right)\left(\vt^2+T_r^2\right)
- \frac\nu\vt \left|\vu\right|^2  -\frac 1 {\sigma\mu^2\vt}\left|\Curl \vc B\right|^2.
\eF

In the sequel, we define $V = \left(\rho,\vu,\vt,E_r,\vc B \right)^T$,
$\overline V = \left(\overline \rho, 0, \overline\vt,\overline{E_r},\overline{\vc B} \right)^T,$ and
\begin{multline}
  \label{eq:def_N(t)}
N(t)^2 = \sup_{0\leq s \leq t} \left\|V(s) - \overline
  V\right\|^2_{H^d\left(\R^3\right)} \\ +\int_0^t \left(\|\Grad V(s)\|^2_{H^{d-1}\left(\R^3\right)}+
\|\Grad \vt(s)\|^2_{H^{d}\left(\R^3\right)} + \|\Grad
E_r(s)\|_{H^{d}\left(\R^3\right)}^2 + \left\|\Grad \vc B(s)\right\|_{H^d(\R^3)}^2 \right)ds \\
+ \int_0^t \left(\left\|\vt(s)
  - T_r(s)\right\|^{2}_{H^{d-1}\left(\R^3\right)} + \left\|\vu(s)\right\|^{2}_{H^{d-1}\left(\R^3\right)}\right)ds.
\end{multline}
Recall that $T_r = E_r^{1/4} a^{-1/4}$. Note also that, since $\Div\left(\vc B\right) = 0$, we have
\begin{equation}\label{eq:10}
  \int_{\R^3}\left|\Curl \vc B\right|^2 = \int_{\R^3} |\Grad \vc B|^2,
\end{equation}
as far as $\vc B\in H^1(\R^3)$, and similarly for any $H^s$ norm. This allows, in the sequel, to replace $\Curl \vc
B$ by $\Grad \vc B$ in all bounds. 
\subsubsection{$L^\infty(H^d)$ estimates}

Using these entropy properties, we are going to prove the following
result:
\begin{Proposition}
  \label{pr:linftyhd_estimate}
Let the assumptions of Theorem~\ref{th:existence_MHD} be
satisfied. Consider a solution $(\vr,\vu,\vt,E_r)$ of system
\eqref{i1bis}-\eqref{i2bis}-\eqref{i3bis}-\eqref{i4bis}-\eqref{i5bis} on $[0,t]$, for some $t>0$.
Then, the energy defined by \eqref{eq:def_N(t)} satisfies
\begin{multline}
  \label{eq:linftyhd_estimate}
 \left\|V(t) -\overline V \right\|_{L^2\left(\R^3\right)}^2 \\+ \int_0^t \left(\left\|\Grad
   \vt(s)\right\|_{L^2\left(\R^3\right)}^2 + \left\|\Grad
   E_r(s)\right\|_{L^2\left(\R^3\right)}^2 + \left\|\vt(s) - T_r(s)\right\|_{L^2\left(\R^3\right)}^2 +
 \left\|\vu(s)\right\|_{L^2\left(\R^3\right)}^2 +  \left\|\Grad \vc B(s)\right\|_{L^2\left(\R^3\right)}^2 
\right)ds  \\ \leq C(N(t)) N(0)^2,
\end{multline}
where the function $C$ is non-decreasing.
\end{Proposition}
\noindent {\bf Proof:} Following the proof of \cite[Lemma 3.1]{LG} we define
\begin{equation}
  \label{eq:def_eta}
  \eta(t,x) = H_{\overline\vt}(\vr,\vt) - \left(\vr - \overline\vr\right)
    \partial_\vr H_{\overline \vt}\left(\overline\vr,\overline\vt\right)
    -H_{\overline \vt}\left(\overline\vr,\overline\vt\right) +
    H_{r,\overline\vt} \left(T_r \right) - H_{r,\overline\vt}\left(\overline T_r\right)..
\end{equation}
We multiply \eqref{ient} by $\overline\vt$, and subtract the result to
\eqref{iener}. Integrating over $[0,t]\times \R^3$, we find
\begin{multline*}
  \int_{\R^3} \left(\frac12 \vr(t) \left|\vu\right|^2(t) +\eta(t,x) + \frac 1 {2\mu} \left|\vc B \right|^2
     \right)dx + \int_0^t \int_{\R^3} \kappa \frac{\overline{\vt}}{\vt^2}
     |\Grad \vt|^2  + \frac{4a}{3\sigma_s} T_r |\Grad
     E_r|^2 \overline{\vt} \\ +\int_0^t \int_{\R^3} \overline\vt \frac{a\sigma_a}{\vt T_r}(\vt+ T_r)\left(\vt
       - T_r\right)^2\left(\vt^2+T_r^2\right) + 
     \frac{\overline \vt}\vt\nu\left|\vu\right|^2 + \dfrac{\overline\vt}{\sigma\mu^2\vt}\left|\Curl \vc B\right|^2 
     \\ 
     \leq \int_{\R^3}\eta(0,x)dx + \int_{\R^3} \vr_0\left|\vec u_0\right|^2 + \frac 1 {2\mu}\int_{\R^3} \left|\vec
       B_0 \right|^2.
\end{multline*}
Defining
\begin{equation}
  \label{eq:M(t)}
  M(t) = \sup_{0\leq s\leq t} \sup_{x\in\R^3}
  \left[\max\left(|\vr(s,x)-\overline\vr|,|\vu(s,x)|,|\vt(s,x)-\overline\vt|,\left|E_r(s,x)-\overline
    {E_r}\right|, \left|\vc B - \overline{\vc B} \right|\right)\right],
\end{equation}
and applying Lemma~\ref{l2}, we find that
\begin{multline*}
    \left\|V(t) -\overline V \right\|_{L^2\left(\R^3\right)}^2 \\+ \int_0^t \left(\left\|\Grad
   \vt(s)\right\|_{L^2\left(\R^3\right)}^2 + \left\|\Grad
   E_r(s)\right\|_{L^2\left(\R^3\right)}^2 + \left\|\vt(s) -
   T_r(s)\right\|_{L^2\left(\R^3\right)}^2 +  \left\|\vu(s)\right\|_{L^2\left(\R^3\right)}^2 + \left\|\Curl \vc B(s)\right\|_{L^2\left(\R^3\right)}^2
\right)ds \\ \leq C(M(t)) N(0),
\end{multline*}
where $C:\R^+\to\R^+$ is non-decreasing. Equation \eqref{eq:10} allows to replace $\Curl \vc B$ by $\Grad \vc B$ in
the above estimate. Finally, we point out that,
since $d>7/2 >3/2$, due to Sobolev embeddings, there exists a universal
constante $C_0$ such that $M(t) \leq C_0 N(t)$. Since $C$ is
non-decreasing, this proves \eqref{eq:linftyhd_estimate}.
\hfill $\square$

\begin{Proposition}
\label{pr:linftyhd-estimates_2}
Setting $V= \left(\vr,\vu,\vt,E_r,\vc B\right)^T$, under the same assumptions as in Theorem~\ref{th:existence_MHD},
we have the following estimate: 
\begin{multline}
  \label{eq:linftyhd_estimate_2}
  \left\|\partial_t V(t) \right\|_{H^{d-1}\left(\R^3\right)} \leq C(N(t))
  \left(\left\|\Grad V\right\|_{H^{d-1}\left(\R^3\right)} + \left\|\Grad
      \vt\right\|_{H^d\left(\R^3\right)} \right. \\ \left.+ \left\|\Grad
      E_r \right\|_{H^d\left(\R^3\right)} + \left\|\vt - T_r\right\|_{H^{d-1}\left(\R^3\right)} +
    \left\|\vu\right\|_{H^{d-1}\left(\R^3\right)} + \left\|\Grad \vc B \right\|_{H^d\left(\R^3\right)}\right).
\end{multline}
\end{Proposition}
\noindent {\bf Proof:} The system satisfied by $V$ may be written formally
\begin{equation}
 \partial_t V + \sum_{j=1}^3 \widehat \cA_j(V) \partial_{x_j} V =
 \widehat  \cD(V) \Delta V - \widehat \cB (V)=0,
\label{systV}
\end{equation}
where
\[
\widehat \cA_1=
\left(
\begin{array}{ccccccccc}
0  &  \vr  &  0  &  0 & 0 & 0 & 0 & 0 & 0\\\
\widehat \alpha' & 0 & 0 & 0 & \widehat \beta' & \widehat \beta'' & 0 & B_2/\vr\mu & B_3/\vr\mu\\\
0 & 0 & 0 & 0 & 0 & 0 & 0 & -B_1/\vr\mu & 0\\\
0 & 0 & 0 & 0 & 0 & 0 & 0 & 0 & -B_1/\vr\mu\\\
0 & \widehat \gamma' & 0 & 0 & 0 & 0 & 0 & 0 & 0\\\
0 & \widehat \gamma'' & 0 & 0 & 0 & u_1 & 0 & 0 & 0\\\
0 & 0 & 0 & 0 & 0 & 0 & u_1 & 0 & 0\\\
0 & B_2 & -B_1 & 0 & 0 & 0 & 0 & u_1 & 0\\\
0 & B_3 & 0 & -B_1 & 0 & 0 & 0 & 0 & u_1
\end{array}
\right),
\]
\medskip
\[
\widehat \cA_2=
\left(
\begin{array}{ccccccccc}
0  & 0 &  \vr  &  0  &  0 & 0 & 0 & 0 & 0 \\\
0 & 0 & 0 & 0 & 0 & 0 & -B_2/\vr\mu & 0 & 0\\\
\widehat \alpha' & 0 & 0 & 0 & \widehat \beta' & \widehat \beta'' & B_1/\vr\mu & 0 & B_3/\vr\mu\\\
0 & 0 & 0 & 0 & 0 & 0 & 0 & 0 & -B_2/\vr\mu\\\
0 & 0 & \widehat \gamma' & 0 & 0 & 0 & 0 & 0 & 0 \\\
0 & 0 & \widehat \gamma'' & 0 & 0 & u_2 & 0 & 0 & 0\\\
0 & -B_2 & B_1 & 0 & 0 & 0 & u_2 & 0 & 0\\\
0 & 0 & 0 & 0 & 0 & 0 & 0 & u_2 & 0\\\
0 & 0 & B_3 & -B_2 & 0 & 0 & 0 & 0 & u_2
\end{array}
\right)
\ \]
\medskip
\[
\widehat \cA_3=
\left(
\begin{array}{ccccccccc}
0 & 0 & 0 & \vr  &  0  &  0 & 0 & 0 & 0\\\
0 & 0 & 0 & 0 & 0 & 0 & -B_3/\vr\mu & 0 & 0\\\
0 & 0 & 0 & 0 & 0 & 0 & 0 & -B_3/\vr\mu & 0\\\
\widehat \alpha' & 0 & 0 & 0 & \widehat \beta' & \widehat \beta'' & B_1/\vr\mu & B_2/\vr\mu & 0\\\
0 & 0 & 0 & \widehat \gamma' & 0 & 0 & 0 & 0 & 0\\\
0 & 0 & 0 & \widehat \gamma'' & 0 & u_3 & 0 & 0 & 0\\\
0 & -B_3 & 0 & B_1 & 0 & 0 & u_3 & 0 & 0\\\
0 & 0 & -B_3 & B_2 & 0 & 0 & 0 & u_3 & 0\\\
0 & 0 & 0 & 0 & 0 & 0 & 0 & 0 & u_3
\end{array}
\right),
\]
and
\[
\widehat \cD=
\left(
\begin{array}{ccccccccc}
0 & 0 & 0 & 0 & 0 & 0 & 0 & 0 & 0\\\
0 & 0 & 0 & 0 & 0 & 0 & 0 & 0 & 0\\\
0 & 0 & 0 & 0 & 0 & 0 & 0 & 0 & 0\\\
0 & 0 & 0 & 0 & 0 & 0 & 0 & 0 & 0\\\
0 & 0 & 0 & 0 & \widehat \delta' & 0 & 0 & 0 & 0\\\
0 & 0 & 0 & 0 & 0 & \widehat \delta'' & 0 & 0 & 0\\\
0 & 0 & 0 & 0 & 0 & 0 &  \lambda & 0 & 0\\\
0 & 0 & 0 & 0 & 0 & 0 & 0 &  \lambda & 0\\\
0 & 0 & 0 & 0 & 0 & 0 & 0 & 0 &  \lambda
\end{array}
\right)
\ \
\widehat \cB:=
\left(
\begin{array}{c}
0 \\\
-\nu\frac{\vu}{\vr}\\\
-\frac{\sigma_a}{\vr C_v}(a\vt^4-E_r)+\frac{1}{\vr C_v}\frac{\lambda}\mu \left|\Curl \vec B\right|^2 + \frac{\nu
  \left|\vec u\right|^2}{\vr C_v}\\\
\sigma_a(a\vt^4-E_r)\\\
0\\\
0 \\\
0
\end{array}
\right),
\]
where
\[
\widehat \alpha'= \frac{p_{\vr}}{\vr},\quad
\widehat \beta'=\frac{p_{\vt}}{\vr},\quad
\widehat \beta''=\frac{1}{3\vr},\quad
\widehat \gamma'=\frac{3\vr p_{\vt}}{3\vr C_v} ,\quad
\widehat \delta'= \frac{\kappa}{\vr C_v},\quad
\widehat \gamma'' = \frac43  E_r,\quad
\widehat \delta''=\frac{1}{3\sigma_s}.
\]
It is possible to symmetrize this nonlinear system in the same spirit as what we have done for the linearized
system (\ref{hrad}). However, we do not need to do so here. So we write
\[
  \partial_t V = -\sum_{j=1}^3 \left[ \widehat \cA_j(V) - \widehat
    \cA_j(\overline V)\right] \partial_{x_j} V - \sum_{j=1}^n \widehat
  \cA_j(\overline V) \partial_{x_j} V + \left[\widehat \cD(V) - \widehat
    \cD(\overline V)\right]\Delta V + \widehat \cD(\overline V) \Delta V
- \widehat \cB (V).
\]
We first observe that these matrices are
Lipschitz continuous with respect to $V$, away from $\vr=0$ and $\vt=0$ and also that
the matrices $\widehat \cB $ and $\widehat \cD$ have, respectively, the
same structure as those defined in \eqref{eq:def_BD}.
Note also that, since $d-1 > 5/2 = 3/2+1,$ Sobolev embeddings imply that
$H^{d-1}\left(\R^3\right)$ is an algebra.
Therefore, we have
\begin{multline*}
  \left\|\partial_t V \right\|_{H^{d-1}\left(\R^3\right)} \leq C_0\left(1+
  \sum_{j=1}^3 \left\|\widehat \cA_j(V) - \widehat
    \cA_j(\overline V)\right\|_{H^{d-1}\left(\R^3\right)}\right)\left\|\Grad V
\right\|_{H^{d-1}\left(\R^3\right)} \\+ C_0\left(1+
  \left\|\widehat \cD(V) - \widehat
   \cD(\overline
    V)\right\|_{H^{d-1}\left(\R^3\right)}\right)\left(\left\|\Delta \vt
\right\|_{H^{d-1}\left(\R^3\right)} +\left\|\Delta E_r\right\|_{H^{d-1}\left(\R^3\right)}+\left\|\Curl\left(\Curl \vc B\right)\right\|_{H^{d-1}\left(\R^3\right)} \right) \\+ C_0\left(1+
  \left\|\widehat \cB (V) - \widehat
   \cB (\overline
    V)\right\|_{H^{d-1}\left(\R^3\right)}\right)\left(\| \vt-T_r\|_{H^{d-1}\left(\R^3\right)} +\| \vu\|_{H^{d-1}\left(\R^3\right)} \right),
\end{multline*}
whence,
\begin{multline*}
  \left\|\partial_t V \right\|_{H^{d-1}\left(\R^3\right)} \leq C_0\left(1+
  \left\|V -\overline  V\right\|_{H^{d-1}\left(\R^3\right)}\right)\left(\left\|\Grad V
\right\|_{H^{d-1}\left(\R^3\right)} +\left\|\Delta \vt
\right\|_{H^{d-1}\left(\R^3\right)} +\left\|\Delta E_r
\right\|_{H^{d-1}\left(\R^3\right)} \vphantom{\left\|\Delta \vc B\right\|_{H^{d-1}\left(\R^3\right)}}\right.\\ \left.+\left\|\Delta \vc B\right\|_{H^{d-1}\left(\R^3\right)}+\| \vt-T_r\|_{H^{d-1}\left(\R^3\right)}+\| \vu\|_{H^{d-1}\left(\R^3\right)}\right),
\end{multline*}
which proves \eqref{eq:linftyhd_estimate_2}.\hfill$\square$
\bigskip

Next, we bound the spatial derivatives as follows:

\begin{Proposition}
  \label{pr:linftyhd-estimates_3} Assume that the hypotheses of
  Theorem~\ref{th:existence_MHD} are satisfied. Let $k\in
  \N^3$ be such that $1\leq |k|\leq d$, where $d>7/2$. Then, we
  have
  \begin{multline}
    \label{eq:linftyhd_estimate_3}
    \left\|\partial_x^k V(t) \right\|_{L^2\left(\R^3\right)}^2 +
    \int_0^t \left(\|\partial_x^k \Grad
\vt(s)\|^2_{L^2\left(\R^3\right)} + \|\partial_x^k \Grad E_r(s)\|_{L^2\left(\R^3\right)}^2 + \left\|\partial_x^k\left(\vt
  - T_r\right)(s)\right\|_{L^2\left(\R^3\right)}^2 \vphantom{+ \left\|\partial_x^k \Grad \vc
B(s)\right\|_{L^2\left(\R^3\right)}^2}\right. \\ \left. + \left\|\partial_x^k \Grad \vc
B(s)\right\|_{L^2\left(\R^3\right)}^2 + \left\|\partial_x^k\vu(s)\right\|_{L^2\left(\R^3\right)}^2 \right)ds \\
\leq C_0 N(0)^2 +C_0 N(t)\int_0^t \left( \left\|\Grad V(s)
\right\|_{H^{d-1}\left(\R^3\right)}^2 + \left\|\Grad
      \vt(s)\right\|_{H^d\left(\R^3\right)}^2 + \left\|\Grad
      E_r \right\|_{H^d\left(\R^3\right)}^2 +\vphantom{+ \left\|\partial_x^k \Grad \vc
B(s)\right\|_{L^2\left(\R^3\right)}^2}\right. \\ \left. \left\|\vt(s) - T_r(s)\right\|_{H^{d-1}\left(\R^3\right)}^2 +
    \left\|\vu(s)\right\|_{H^{d-1}\left(\R^3\right)}^2 + \left\|\Grad \vc B(s) \right\|_{H^d\left(\R^3\right)}^2\right)ds
  \end{multline}
\end{Proposition}
\noindent {\bf Proof:} Here, we need to symmetrize the nonlinear system. For this purpose, we multiply
(\ref{systV}) on the right by the matrix 
\begin{equation}
\widehat \cA_0(V)=
\left(
\begin{array}{ccccccccc}
\mu\frac{\widehat \alpha'}{\vr} & 0 & 0 & 0 & 0 & 0 & 0 & 0 & 0\\\
0 & \mu & 0 & 0 & 0 & 0 & 0 & 0 & 0\\\
0 & 0 & \mu & 0 & 0 & 0 & 0 & 0 & 0\\\
0 & 0 & 0 & \mu & 0 & 0 & 0 & 0 & 0\\\
0 & 0 & 0 & 0 & \mu\frac{\widehat \beta'}{\widehat \gamma'} & 0 & 0 & 0 & 0\\\
0 & 0 & 0 & 0 & 0 & \mu\frac{\widehat \beta''}{\widehat \gamma''} & 0 & 0 & 0\\\
0 & 0 & 0 & 0 & 0 & 0 & 1 & 0 & 0\\\
0 & 0 & 0 & 0 & 0 & 0 & 0 & 1 & 0\\\
0 & 0 & 0 & 0 & 0 & 0 & 0 & 0 & 1
\end{array}
\right).
\label{A0}
\end{equation}
This gives
\begin{equation}
  \label{eq:13}
  \widehat \cA_0(V) \partial_t V = -\sum_{j=1}^3  \widecheck \cA_j(V)  \partial_{x_j} V  + \widecheck \cD(V) \Delta V 
- \widecheck \cB (V)= 0,
\end{equation}
where $\widecheck \cA_j(V) = \widehat \cA_0(V) \widehat \cA_j(V)$, $\widecheck
\cB (V) = \widehat\cA_0(V) \widehat\cB (V),$ and $\widecheck \cD(V) =
\widehat \cA_0(V) \widehat \cD(V)$ are all symmetric matrices. Applying $\partial_x^k$ to (\ref{eq:13}) then taking the scalar
product with the vector $\partial_x^k V$, and integrating over
$[0,t]\times \R^3$, we find
\begin{multline*}
  \frac12 \int_{\R^3} \left[\left( \widehat \cA_0 \partial_x^k V\right)
    \cdot \partial_x^k V \right]^t_0dx + \int_0^t \int_{\R^3}
  \left(\widecheck\cD\Grad \left(\partial_x^k V\right)\right)\cdot\Grad
  \left(\partial_x^k V\right)dxdt + \int_0^t \int_{\R^3} \left(
    \widecheck\cB(V)\partial_x^k V \right)\cdot
  \left(\partial_x^k V \right)dxdt \\
= \int_0^t \int_{\R^3} \left( \frac12\left(I_1+I_2\right) - I_3 - I_4 - I_5\right)dxdt,
\end{multline*}
where
  \begin{align*}
  I_1 &= \partial_t\left(\widehat\cA_0(V)\right)\partial_x^k V
  \cdot \partial_x^k V, & I_2 &= \sum_{j=1}^3 \partial_{x_j}
  \left(\widecheck\cA_j(V)\right) \partial_x^k V \cdot \partial_x^k V,
  & I_3 &= \left[\partial_x^k , \widehat\cA_0(V)\right] \partial_t V
  \cdot \partial_x^k V, \\ I_4 &= \sum_{j=1}^3
  \left[\partial_x^k , \widecheck\cA_j(V)\right] \partial_{x_j} V
  \cdot\partial_x^k V, & I_5&
  = \partial_x^k\left(\widecheck\cB (V)\right)  \cdot \partial_x^k V.
  \end{align*}
We estimate separately each term of the right-hand side.

 First, we have
\[
  \int_0^t \int_{\R^3}  |I_1| \leq C\int_0^t \int_{\R^3}
  \left|\partial_x^k V\right|^2 \left|\partial_t V \right| 
\leq  C\int_0^t \int_{\R^3}
  \left|\partial_x^k V\right|^2 \left(\left|\Grad V \right| +\left|\cB (V)\right| + \left|\cD\Delta V\right|\right)
\]
\[
\leq  C N(t) \int_0^t \left\|\partial_x^k
  V(s)\right\|^2_{L^2\left(\R^3\right)} ds,
\]
where we have used Sobolev embeddings and the fact that $d>7/2$. A
similar computation gives
\begin{displaymath}
  \int_0^t \int_{\R^3}  |I_2| \leq C N(t) \int_0^t \left\|\partial_x^k
  V(s)\right\|^2_{L^2\left(\R^3\right)} ds.
\end{displaymath}
We estimate $I_3$ by applying Cauchy-Schwarz inequality:
\begin{displaymath}
  \int_0^t \int_{\R^3} \left|I_3\right| \leq \int_0^t \left\|\partial_x^k
  V\right\|_{L^2\left(\R^3\right)} \left\| \left[\partial_x^k,\widehat \cA_0(V)\right] \partial_t
  V\right\|_{L^2\left(\R^3\right)}.
\end{displaymath}
Then, we apply the same estimate for commutators and composition of
functions (see \cite[Proposition 2.1]{M}), and $|k|\leq d$:
\begin{multline*}
  \left\| \left[\partial_x^k,\widehat \cA_0(V)\right] \partial_t
  V\right\|_{L^2\left(\R^3\right)} = \left\|
  \left[\partial_x^k,\widehat\cA_0(V)-\widehat \cA_0(\overline V)\right] \partial_t
  V\right\|_{L^2\left(\R^3\right)}\\ \leq C \left(\left\|\partial_t V
  \right\|_{L^\infty\left(\R^3\right)} \left\|\Grad \widehat\cA_0(V)\right\|_{H^{d-1}\left(\R^3\right)}
  + \left\|\partial_t V
  \right\|_{H^{d-1}\left(\R^3\right)} \left\|\Grad \widehat\cA_0(V)\right\|_{L^\infty\left(\R^3\right)}\right).
\end{multline*}
Moreover, we have
\begin{displaymath}
  \left\|\Grad \widehat\cA_0(V)\right\|_{H^{d-1}\left(\R^3\right)} \leq C \left\| V -
    \overline V \right\|_{H^d\left(\R^3\right)} \leq C N(t),
\end{displaymath}
and
\begin{displaymath}
  \left\|\Grad \widehat\cA_0(V)\right\|_{L^\infty\left(\R^3\right)} \leq C \left\| \Grad V \right\|_{H^{d-1}\left(\R^3\right)} \leq C N(t).
\end{displaymath}
Hence, $I_3$ satisfies
\begin{multline*}
  \int_0^t \int_{\R^3} \left|I_3\right| \leq C(N(t)) N(t) \int_0^t
  \left(\|\Grad V(s)\|_{H^{d-1}\left(\R^3\right)}^2 +\|\Grad \vt(s)\|_{H^d\left(\R^3\right)}^2+ \left\|\Grad
      E_r(s) \right\|_{H^d\left(\R^3\right)}^2 \vphantom{\left\|\Grad \vc B \right\|_{H^d\left(\R^3\right)}}
  \right. \\ \left.
+ \left\|\vt(s) - T_r(s)\right\|_{H^{d-1}\left(\R^3\right)}^2 +
    \left\|\vu(s)\right\|_{H^{d-1}\left(\R^3\right)}^2 + \left\|\Grad \vc B(s) \right\|_{H^d\left(\R^3\right)}^2\right)ds.
\end{multline*}
Here, we have used \eqref{eq:linftyhd_estimate_2}.

 The integral of $I_4$
is dealt with using similar computations.

 Turning to $I_5$, we use the
particular form of $\partial_x^k \widecheck\cB (V)$. More precisely, we have
\[
  \partial_x^k \left(\widecheck\cB (V)\right)  \cdot\partial_x^k V =
\partial_x^k \left(\frac{\vec f_m}{\vr}-\nu\frac{\vu}{\vr}\right) \cdot\partial_x^k \vu
-\partial_x^k \left(\frac{\sigma_a}{\vr}(a\vt^4-E_r)+\frac{E_m-\vec f_m\cdot\vu}{\vr}\right) \cdot\partial_x^k \vt
\]
\[
+\partial_x^k \left(\sigma_a(a\vt^4-E_r)\right) \cdot\partial_x^k E_r,
\]
from which, using estimates for composition of functions (see Proposition 2.1 in \cite{M}) we infer
\begin{displaymath}
  \int_0^t \int_{\R^3} \left|I_5\right| \leq C N(t)\int_0^t \left\|\partial_x^k
  V(s)\right\|_{L^2\left(\R^3\right)}^2 ds.
\end{displaymath}
Collecting the estimates on $I_1$, $I_2$, $I_3$, $I_4$ and $I_5$, we have
proved \eqref{eq:linftyhd_estimate_3}.\hfill$\square$

\medskip

The above results allow to derive the following bound:
\begin{Proposition}\label{pr:fin_linftyhd_estimate}
  Assume that the assumptions of Theorem~\ref{th:existence_MHD}
  are satisfied. Then, there exists a non-decreasing function $C:\R^+\to
  \R^+$ such that
  \begin{multline}
    \label{eq:fin_linfty_hd}
    \left\|V - \overline V \right\|_{H^d\left(\R^3\right)}^2 + \int_0^t\left(\left\|\Grad
   \vt(s)\right\|_{H^d\left(\R^3\right)}^2 + \left\|\Grad
   E_r(s)\right\|_{H^d\left(\R^3\right)}^2 + \left\|\vt(s) - T_r(s)\right\|_{H^d\left(\R^3\right)}^2\vphantom{+ \left\|\Grad \vc B \right\|_{H^d\left(\R^3\right)}^2} \right. \\ \left.+
 \left\|\vu(s)\right\|_{H^d\left(\R^3\right)}^2 + \left\|\Grad \vc B(s) \right\|_{H^d\left(\R^3\right)}^2
\right)ds \\
    \leq C(N(t)) \left[N(0)^2 + N(t) \int_0^t \left( \left\|\Grad V(s)
\right\|_{H^{d-1}\left(\R^3\right)}^2 + \left\|\Grad
      \vt(s)\right\|_{H^d\left(\R^3\right)}^2 + \left\|\Grad
      E_r(s) \right\|_{H^d\left(\R^3\right)}^2 +\vphantom{+ \left\|\partial_x^k \Grad \vc
B(s)\right\|_{L^2\left(\R^3\right)}^2}\right.\right. \\ \left.\left. \left\|\vt(s) - T_r(s)\right\|_{H^{d-1}\left(\R^3\right)}^2 +
    \left\|\vu(s)\right\|_{H^{d-1}\left(\R^3\right)}^2 + \left\|\Grad \vc B(s) \right\|_{H^d\left(\R^3\right)}^2\right)ds \right].
  \end{multline}
\end{Proposition}
\noindent {\bf Proof:} We sum up estimates \eqref{eq:linftyhd_estimate_3} over all
multi-indices $k$ such that $|k|\leq d$, and add this to
\eqref{eq:linftyhd_estimate}. This leads to
\eqref{eq:fin_linfty_hd}. \hfill $\square$

\subsubsection{$L^2(H^{d-1})$ estimates}

In this section, we derive bounds on the right-hand side of
\eqref{eq:fin_linfty_hd}. For this purpose, we adapt the strategy of
\cite{SK}, which was further developed in \cite{HN}. We apply the
Fourier transform to the linearized system and use the compensating
matrix $K$ to prove estimates on the space derivatives of $V$.

\begin{Proposition}
  \label{pr:l2-estimates} Assume that the assumptions of
  Theorem~\ref{th:existence_MHD} are satisfied. Then there
  exists a non-decreasing function $C:\R^+ \to \R^+$ such that
  \begin{equation}
    \label{eq:l2-estimate}
    \int_0^t \left\|\Grad V (s)\right\|_{H^{d-1}\left(\R^3\right)} ds
    \leq C(N(t)) \left(N(t) + \left\|V_0 - \overline V\right\|_{H^d\left(\R^3\right)} \right)
  \end{equation}
\end{Proposition}
\noindent {\bf Proof:} As a first step, we apply the symmetrizer of the
linearized system \eqref{hrad} (which leads to \eqref{hsym}) to the
nonlinear system \eqref{i1bis}-\eqref{i2bis}-\eqref{i3bis}, which then
reads
\begin{displaymath}
  \widetilde \cA_0(V)\partial_t V + \sum_{j=1}^3 \widetilde
  \cA_j(V) \partial_{x_j} V  = \widetilde \cD\Delta V - \widetilde \cB (V) V.
\end{displaymath}
Of course, this system is not symmetric. However, the corresponding
linearized system \eqref{hsym} is symmetric. 
Next, we rewrite the nonlinear system by setting $U = V-\overline V$:
\begin{displaymath}
 \widetilde \cA_0(V)\partial_t U + \sum_{j=1}^3 \widetilde
  \cA_j(V) \partial_{x_j} U  = \widetilde \cD\Delta U - \widetilde \cB (V) U -
  \widetilde \cB (V)\overline V.
\end{displaymath}
Therefore, multiplying this system on the left by $\widetilde \cA_0(\overline V)
\left(\widetilde \cA_0(V)\right)^{-1},$ we find
\begin{equation}\label{eq:2}
  \widetilde \cA_0(\overline V) \partial_t U + \sum_{j=1}^3
  \widetilde \cA_j(\overline V) \partial_{x_j} U = H,
\end{equation}
where
\begin{multline*}
  H = - \widetilde \cA_0(\overline V) \sum_{j=1}^3\left[\left(\widetilde
      \cA_0(V)\right)^{-1} \widetilde \cA_j(V) -
      \left(\widetilde \cA_0(\overline V)\right)^{-1}\widetilde \cA_j(\overline V)
    \right]\partial_{x_j} V \\
+\widetilde \cA_0(\overline V)\left(\widetilde \cA_0(V)\right)^{-1}
\widetilde \cD\Delta U - \widetilde \cA_0(\overline V)\left(\widetilde \cA_0(V)\right)^{-1}
\widetilde \cB (V) U - \widetilde \cA_0(\overline V)\left(\widetilde
  \cA_0(V)\right)^{-1} \widetilde \cB (V) \overline V.
\end{multline*}
We apply the Fourier transform to \eqref{eq:2}, and then multiply on the
left by $-i\left(\widehat U\right)^*K\left(\frac{\xi}{|\xi|}\right)$,
where $~^*$ denotes the transpose of the complex conjugate, and $K$ is the compensating
matrix (see Proposition~\ref{pr:compensating_matrix}). Taking the real
part of the result, we infer
\begin{equation}
  \label{eq:1}
\im \left(\left(\widehat U\right)^* K\left(\frac{\xi}{|\xi|}\right)
  \cA_0\left(\overline V \right) \partial_t \widehat U\right) +
|\xi|\left(\widehat U\right)^*
K\left(\frac{\xi}{|\xi|}\right)\cA\left(\frac{\xi}{|\xi|}\right) \widehat
U = \im \left( \left(\widehat U\right)^* K\left(\frac{\xi}{|\xi|}\right)
\widehat H\right),
\end{equation}
where the matrix $\cA\left(\frac{\xi}{|\xi|}\right)$ is defined by
\eqref{eq:3}. According to Proposition~\ref{pr:compensating_matrix},
$K\cA_0(\overline V)$ is skew-symmetric, hence
\begin{displaymath}
  \im \left(\left(\widehat U\right)^* K\left(\frac{\xi}{|\xi|}\right)
  \cA_0\left(\overline V \right) \partial_t \widehat U\right) = \frac12
\frac{d}{dt} \im \left( \left(\widehat U\right)^* K\left(\frac{\xi}{|\xi|}\right)
  \cA_0\left(\overline V \right)  \widehat U\right).
\end{displaymath}
Next, we also have
\begin{equation}
  \label{eq:4}
|\xi| \left(\widehat U\right)^* K\left(\frac{\xi}{|\xi|}\right)
  \cA\left(\frac{\xi}{|\xi|} \right)  \widehat U = |\xi|\left(\widehat
    U\right)^*\left[
 K\left(\frac{\xi}{|\xi|}\right)
 \cA\left(\frac{\xi}{|\xi|} \right) + \cB \left(\frac\xi{|\xi|}\right)
\right]  \widehat U - |\xi|\left(\widehat U\right)^* \widetilde \cB 
\widehat U - |\xi| \left(\widehat U\right)^* \widetilde \cD
\widehat U.
\end{equation}
Hence, still applying Proposition~\ref{pr:compensating_matrix}, there
exists $\alpha_1>0$ and $\alpha_2>0$ such that
\begin{multline}
  \label{eq:5}
|\xi| \left(\widehat U\right)^* K\left(\frac{\xi}{|\xi|}\right)
  \cA\left(\frac{\xi}{|\xi|} \right)  \widehat U \geq \alpha_1
  |\xi|\left|\widehat U\right|^2 \\- \alpha_2 \frac 1 {|\xi|}
  \left(\left|\xi\left(\widehat{\vt-\overline\vt}\right)\right|^2 +
    \left|\xi \left(\widehat{E_r-\overline{E_r}}\right)\right|^2 + 
    |\xi|^2\left| \widehat{\vc B - \overline{\vc B}}\right|^2 + 
     |\xi|^2\left|\widehat{\vu}\right|^2  + |\xi|^2\left|\widehat{\vt - T_r}\right|^2\right).
\end{multline}
Finally, we estimate the right-hand side of \eqref{eq:1} using
Cauchy-Schwarz inequality and Young inequality:
\begin{equation}
  \label{eq:6}
\left| \im \left( \left(\widehat U\right)^* K\left(\frac{\xi}{|\xi|}\right)
\widehat H\right)\right| \leq \varepsilon |\xi|\left|\widehat U\right|^2
+ C_\varepsilon \frac 1 {|\xi|} \left|\widehat H\right|^2,
\end{equation}
for any $\varepsilon>0$. We choose $\varepsilon$ small enough, insert
\eqref{eq:4}-\eqref{eq:5}-\eqref{eq:6} into \eqref{eq:1}, and find
\begin{multline*}
  |\xi| \left|\widehat U\right|^2 \leq C \left[\frac 1 {|\xi|}
  \left(\left|\xi \left(\widehat{\vt-\overline\vt}\right)\right|^2 +
    \left|\xi \left(\widehat{E_r-\overline{E_r}}\right)\right|^2 +
    |\xi|^2\left|\widehat{\vc B - \overline{\vc B}}\right|^2+|\xi|^2\left|\widehat{\vu}\right|^2 + |\xi|^2\left|\widehat{\vt - T_r}\right|^2\right) \right.\\ \left.+\frac 1 {|\xi|}
  \left|\widehat H\right|^2 - \frac{d}{dt} \im \left( \left(\widehat U\right)^* K\left(\frac{\xi}{|\xi|}\right)
  \cA_0\left(\overline V \right)  \widehat U\right)  \right].
\end{multline*}
We multiply this inequality by $|\xi|^{2l-1}$, for some $1\leq l\leq d$,
and get
\begin{multline}\label{eq:7}
  |\xi|^{2l} \left|\widehat U\right|^2 \leq C \left[|\xi|^{2l-2}
  \left(\left|\xi \left(\widehat{\vt-\overline\vt}\right)\right|^2 +
    \left|\xi \left(\widehat{E_r-\overline{E_r}}\right)\right|^2 +|\xi|^2\left|\widehat{\vc B - \overline{\vc B}}\right|^2+
    |\xi|^2\left|\widehat{\vu}\right|^2 + |\xi|^2\left|\widehat{\vt - T_r}\right|^2\right) +|\xi|^{2l-2}
  \left|\widehat H\right|^2 \right. \\ \left.- |\xi|^{2l-1}\frac{d}{dt} \im \left( \left(\widehat U\right)^* K\left(\frac{\xi}{|\xi|}\right)
  \cA_0\left(\overline V \right)  \widehat U\right)  \right].
\end{multline}
We integrate this inequality over $[0,t]\times \R^3$, and use
Plancherel's theorem:
\begin{multline}\label{eq:8}
  \int_0^t \int_{\R^3} \sum_{|k|=l-1} \left|\partial_x^k \Grad V
  \right|^2 \leq C\int_0^t \int_{\R^3}
  \sum_{|k|=l-1}\left(\left|\partial_x^k \Grad \vt\right|^2 +
    \left|\partial_x^k \Grad E_r\right|^2  + \left|\partial_x^k \Grad\vc B\right|^2 +
    \left|\partial_x^k \Grad\vu\right|^2 + \left|\partial_x^k H\right|^2\right) \\
 + C \im \int_{\R^3} |\xi|^{2l-1}\left[ \left(\widehat U\right)^* K\left(\frac{\xi}{|\xi|}\right)
  \cA_0\left(\overline V \right)  \widehat U\right]^t_0.
\end{multline}
The matrix $K\left(\frac{\xi}{|\xi|}\right)$ is uniformly bounded for
$\xi\in \R^3\setminus\{0\}$, so we have
\begin{multline*}
  \im \int_{\R^3} |\xi|^{2l-1}\left[ \left(\widehat U\right)^* K\left(\frac{\xi}{|\xi|}\right)
  \cA_0\left(\overline V \right)  \widehat U\right]^t_0
\leq C \left(\int_{\R^3} \left(1+|\xi|^2\right)^l \left|\widehat
    U(t)\right|^2 + \int_{\R^3} \left(1+|\xi|^2\right)^l \left|\widehat
    U_0\right|^2 \right) \\ \leq C\left(\left\|V - \overline
    V\right\|_{H^l\left(\R^3\right)}^2 + \left\|V_0 - \overline V\right\|_{H^l\left(\R^3\right)}^2 \right)
\end{multline*}
We insert this estimate into \eqref{eq:8}, sum the result over $1\leq l\leq d$, which leads to
\begin{multline}\label{eq:9}
  \int_0^t \left\|\Grad V\right\|_{H^{d-1}\left(\R^3\right)} \leq C
  \left( \left\|V - \overline
    V\right\|_{H^d\left(\R^3\right)}^2 + \left\|V_0 - \overline V\right\|_{H^d\left(\R^3\right)}^2
  \right. \\ \left.+ \int_0^t \left(\left\|\Grad \vt\right\|^2_{H^{d-1}\left(\R^3\right)} +
    \left\|\Grad E_r\right\|^2_{H^{d-1}\left(\R^3\right)} + \left\|\Grad
      \vu\right\|^2_{H^{d-1}\left(\R^3\right)}  + \left\|\Grad
      \vc B\right\|^2_{H^{d-1}\left(\R^3\right)}  + \|H\|_{H^{d-1}\left(\R^3\right)}^2 \right)\right).
\end{multline}
In order to conclude, we need to estimate the perturbation $H$. For this
purpose, we use that $H^{d-1}\left(\R^3\right)$ is an algebra: for any $s\leq t$,
\begin{displaymath}
  \|H(s)\|_{H^{d-1}\left(\R^3\right)}^2 \leq C N(t) \left\|\Grad
    V\right\|_{H^{d-1}\left(\R^3\right)}.
\end{displaymath}
Inserting this into \eqref{eq:9}, we prove \eqref{eq:l2-estimate}.\hfill $\square$

\medskip

We are now in position to conclude with the


{\bf Proof of Theorem~\ref{th:existence_MHD}:} We first point
out that local existence for system
\eqref{i1bis}-\eqref{i2bis}-\eqref{i3bis} may be proved using standard
fix-point methods. We refer to \cite{M} for the proof. The existence is
proved in the following functional space:
\begin{multline*}
  X(0,T) = \left\{ V , \quad V-\overline V \in
    C\left([0,T];H^d\left(\R^3\right)\right), \quad \Grad V \in L^2\left([0,T];
      H^{d-1}\left(\R^3\right)\right), \right. \\\left. \Grad \vt, \Grad E_r, \Grad \vc B \in
    L^2\left([0,T]; H^d\left(\R^3\right)\right)\right\}.
\end{multline*}
In order to prove global existence, we argue by contradiction, and
assume that $T_c>0$ is the maximum time existence. Then, we necessarily
have
\begin{displaymath}
  \lim_{t\to T_c} N(t) = +\infty,
\end{displaymath}
where $N(t)$ is defined by \eqref{eq:def_N(t)}. We are thus reduced to
prove that $N$ is bounded. For this purpose, we use the method of
\cite{LG}, which was also used in \cite{nishida}. First note  that, due to
Proposition~\ref{pr:fin_linftyhd_estimate} on the one hand, and to
Proposition~\ref{pr:l2-estimates} on the other hand, we know that there
exists a non-decreasing continuous function $C:\R^+ \to \R^+$ such that
\begin{equation}\label{eq:12}
\forall T\in [0,T_c], \quad  N(t)^2 \leq C(N(t))\left(N(0)^2 + N(t)^3\right).
\end{equation}
Hence, setting $N(0) = \varepsilon$, we have
\begin{equation}\label{eq:11}
  \frac{N(t)^2}{\varepsilon^2 + N(t)^3} \leq C(N(t)),
\end{equation}
Studying the variation of $\phi(N) = N^2/\left(\varepsilon^2 +
  N^3\right)$, we see 
 that $\phi'(0) = 0$, that $\phi$ is
increasing on the interval $\left[0, \left(2\varepsilon^2\right)^{1/3}\right]$ and decreasing on the
interval $\left[\left(2\varepsilon^2\right)^{1/3},+\infty
\right)$. Hence,
\begin{displaymath}
  \max \phi = \phi\left(\left(2\varepsilon^2\right)^{1/3} \right) =
  \frac 13 \left(\frac 2 \varepsilon\right)^{2/3}.
\end{displaymath}
Hence, the function $C$ being independent of $\varepsilon$, we can
choose $\varepsilon$ small enough to have $\phi(N)\leq C(N)$ for all
$N\in [0,N^*]$, where $N^*>0$. Since $C$ is continuous, \eqref{eq:11}
implies that $N\leq N^*$. This is clearly in contradiction with
\eqref{eq:12}.\hfill $\square$

\vskip0.25cm
{\bf Acknowlegment:}
{\it \v S\' arka Ne\v casov\' a acknowledges the support of the GA\v CR (Czech Science Foundation) project  17-01747S  in the framework of RVO: 67985840.

 Bernard Ducomet is partially supported by the ANR project INFAMIE (ANR-15-CE40-0011).}


\vskip0.25cm
\centerline{Xavier Blanc}
 \centerline{Univ. Paris Diderot, Sorbonne Paris Cit\'e,}
\centerline{Laboratoire Jacques-Louis Lions,}
\centerline{UMR 7598, UPMC, CNRS, F-75205 Paris, France}
 \centerline{E-mail: blanc@ann.jussieu.fr}
\vskip0.25cm
\centerline{Bernard Ducomet}
 \centerline{CEA, DAM, DIF}
\centerline{ F-91297 Arpajon, France}
 \centerline{E-mail: bernard.ducomet@cea.fr}
\vskip0.25cm
\centerline{\v S\' arka Ne\v casov\' a}
\centerline{Institute of Mathematics of the Academy of Sciences of the Czech Republic}
\centerline{\v Zitn\' a 25, 115 67 Praha 1, Czech Republic}
\centerline{E-mail: matus@math.cas.cz}

    \end{document}